\documentclass[12pt,a4paper]{amsart}
\usepackage{latexsym}
\usepackage{amssymb}

\usepackage{amssymb, verbatim}
\usepackage{graphics}  % for the diagrams!
\usepackage{graphpap}  % for the diagrams!
\usepackage{amssymb, epic, xypic} % for the diagrams!

\makeatletter
\renewcommand{\@begintheorem}[2]{
\rm \trivlist \item [\hskip \labelsep {\bf #2\ \ #1.}]
                                }
\makeatother

\makeatletter

\DeclareFontFamily{U}{cyr}{}
\DeclareFontShape{U}{cyr}{m}{n}{
  <5> wncyr5 <6> wncyr6 <7> wncyr7 <8> wncyr8 <9> wncyr9 <10->
wncyr10}{}
\DeclareMathAlphabet{\mathcyr}{U}{cyr}{m}{n}

\input cyracc.def
\input epsf

\newcommand{\ww}{{\omega}}

\newcommand{\ts}{\vspace{\baselineskip}\noindent{\bf Proof.}$\;\;$}
\newcommand{\ZZ}{{\bf Z}}
\newcommand{\QQ}{{\bf Q}}
\newcommand{\RR}{{\bf R}}
\newcommand{\CC}{{\bf C}}

\newcommand{\FF}{{\bf F}}

\newcommand{\PP}{{\bf P}}

\newcommand{\ccA}{{\mathcal A}}

\newcommand{\ccH}{{\mathcal H}}

\newcommand{\ccS}{{\mathcal H}}

\newcommand{\ccX}{{\mathcal X}}

\newcommand{\bes}{\begin{equation*}}
\newcommand{\ees}{\end{equation*}}

\newcommand{\thchar}[8]{\ensuremath{ \left[\begin{array}{cccc}
\negthickspace#1\negthickspace&\negthickspace#2\negthickspace&\negthickspace#3\negthickspace&\negthickspace#4\negthickspace\\\negthickspace#5\negthickspace&\negthickspace#6\negthickspace&\negthickspace#7\negthickspace&\negthickspace#8\negthickspace\end{array}\right]}}

\def\DynkinEEEE#1#2#3#4#5#6#7
   {\begin{picture}(150, 55)%
   \put(5, 40){\circle{4}}%
   \put(7, 40){\line(1, 0){28}}%
   \put(37, 40){\circle{4}}%
   \put(39, 40){\line(1, 0){28}}%
   \put(69, 40){\circle{4}}%
   \put(71, 40){\line(1, 0){28}}%
   \put(101, 40){\circle{4}}%
    \put(103, 40){\line (1,0){28}}%
    \put(133, 40){\circle{4}}%
   \put(69, 38){\line(0, -1){18}}%
   \put(69, 18){\circle{4}}%
   \put(5, 46){\makebox(0, 0)[b]{$#2$}}%
   \put(37, 46){\makebox(0, 0)[b]{$#3$}}%
   \put(69, 46){\makebox(0, 0)[b]{$#4$}}%
   \put(101, 46){\makebox(0, 0)[b]{$#5$}}%
   \put(133, 46){\makebox(0, 0)[b]{$#6$}}%
   \put(69, 13){\makebox(0, 0)[t]{$#1$}}%
\put(250, 36){\makebox(0, 0)[lb]{$#7$}}%
\end{picture}}

\title{A Picard modular fourfold and the Weyl group $W(E_6)$}
\author{Bert van Geemen}
\author{Kenji Koike}
\address{Dipartimento di Matematica, Universit\`a di Milano,
Via Saldini 50, 20133 Milano, Italia}
\address{Faculty of Education, University of Yamanashi, Takeda 4-4-37, 
Kofu, Yamanashi 400-8510, Japan}
\begin{document}

\begin{abstract}
We study the geometry of a Picard modular fourfold which parametrizes
abelian fourfolds of Weil type for the field of cube roots of unity.
We find a projective model of this fourfold as a singular, degree ten, hypersurface $\ccX$ in projective 5-space. 
The 
Weyl group $W(E_6)$ acts on $\ccX$ and we provide an explicit description of this action. Moreover, we describe various special subvarieties as well as the boundary of $\ccX$.
\end{abstract}

\maketitle

\section*{Introduction}
The Picard modular fourfold which we consider parametrizes principally polarized abelian varieties of dimension four with an automorphism of order three. 
The period matrices of these abelian varieties are the fixed points of an element $M$ of order three in $Sp(8,\ZZ)$. 
Using second order thetanulls, we show this fixed point locus is mapped to a projective space of dimension five. To find the equation for the image, we use classical relations between thetanulls (and a computer!). It turns out that the closure of the image, denoted by $\ccX$, is a hypersurface of degree ten. 

The elements in $Sp(8,\ZZ)$ which normalize the subgroup generated by $M$
act by projective transformations on $\ccX$. We show that they generate 
a subgroup of $Aut(\ccX)$, which is isomorphic to the Weyl group $W(E_6)$ of the root system $E_6$. Actually the action of $W(E_6)$ on the $\PP^5$ is induced from its standard representation. 

After having established these basic facts, 
we consider fixed point loci in $\ccX$ of elements of $W(E_6)$. 
These correspond to abelian fourfolds with further automorphisms. 
The main  problem we encountered was to find explicit elements in 
$Sp(8,\ZZ)$ which represent these automorphisms. Equivalently, we had to find elements, of finite order and up to conjugacy, in the normalizer of $M$ which map to a given subgroup in $W(E_6)$.
There seems to be no systematic way to proceed, but in all the examples we succeeded in finding them.

As a result of our efforts, we found that the Hessian 
of Igusa's quartic is a Shimura variety. In fact it is the quotient of the Siegel upper half space $\ccS_2$ by an arithmetic group.
This was conjectured by Bruce Hunt. We also find a projective model of a moduli space of abelian fourfolds whose endomorphism algebra contains the field of $12$-th roots of unity.

The methods we used are based on those from \cite{pmv}. 
The computations involved in this case are however rather more demanding and we used the computer algebra system Magma \cite{magma} extensively.

More intrinsically, the Picard modular fourfold we study is the moduli space of abelian fourfolds of Weil type for the field of cube roots of unity. In fact, $M$ induces an automorphism of order three which has two distinct eigenvalues, each with multiplicity two, on the tangent space at the origin of these abelian fourfolds. Such  abelian fourfolds are interesting since they have exceptional Hodge classes, shown to be algebraic in this case by C.\ Schoen \cite{chad}. Moreover, the second cohomology group of these abelian fourfolds has a Hodge substructure of K3-type (see \cite{lom}).
It is not yet known if the Kuga-Satake 
correspondence is realized by an algebraic cycle. 
The singular variety $\ccX$ is a projective model of a Shimura variety 
associated to the group $SU(2,2)$. It is an interesting problem to find the Hodge numbers $h^{p,0}$, which are birational invariants, for a(ny) desingularization of $\ccX$. The regular four-forms (if any) on such a desingularization would provide examples of holomorphic modular forms on an arithmetic subgroup of $SU(2,2)$.
We hope to return to these matters in the future.

\section{Abelian varieties and theta functions}

\subsection{The automorphism}\label{auto}
The group $Sp(8,\ZZ)$ 
(also written as $\Gamma_4$) 
of $8\times 8$ symplectic matrices with integer coefficients,
is defined as
$$
Sp(8,\ZZ)\,:=\,\{\,M\,\in\,M_8(\ZZ)\,:\quad ME\,{}^t\!M\,=\,E\,\}~,
\quad E\,:=\,\left(\begin{array}{cc}0& I\\-I&0 \end{array}\right)~,
$$
where $I$ is the $4\times4$ identity matrix.
For $\tau\in\ccS_4$, the Siegel space of $4\times 4$ complex symmetric matrices with 
positive definite imaginary part, 
one defines the principally polarized abelian variety (ppav)
$$
X_\tau\,\cong\,\CC^4/(\ZZ^8\Omega_\tau),\qquad \Omega_\tau\,:=
\,\left(
\begin{array}{cc} \tau\\ I \end{array}\right),\qquad(\tau\in \ccS_4)~,
$$
so we consider the elements of $\ZZ^8$, $\CC^4$ as row vectors.
The symplectic group $Sp(8,\ZZ)$ acts on $\ccS_4$, 
the action of a matrix $N$ with blocks $a,b,c,d$ on $\ccS_4$ is given by 
$N\cdot\tau:=(a\tau+b)(c\tau+d)^{-1}$ as usual.

To define the abelian fourfolds with an automorphism of order three we introduce the matrices:
$$
A\,:=\,\left(\begin{array}{cc}-I& -I\\I&0 \end{array}\right)\quad(\in GL(4,\ZZ)),
\qquad
M\,:=\,\,\left(\begin{array}{cccc}A&0\\0&{}^tA^{-1} \end{array}\right)
\quad(\in Sp(8,\ZZ))~,
$$
where now $I$ is the $2\times 2$ identity matrix. 
Both $A$ and $M$ satisfy the equation $x^2+x+1=0$. 
The fixed point locus of the matrix $M$ above is denoted by
$$
\ccS_4^M\,=\,\{\tau\,\in\,\ccS_4\,:\,M\cdot \tau\,=\,\tau\,\}
\,=\,
\left\{\left(\begin{array}{cc}b+{}^tb& -b\\-{}^tb& b+{}^tb 
\end{array}\right)\,\in\,\ccS_4\,\right\}~,
$$
where $b$ is a $2\times 2$ complex matrix, 
and we used that $M\cdot\tau=\tau$ is equivalent to $A\tau{}^t\!A=\tau$, 
which again is equivalent to $A\tau=\tau({}^tA)^{-1}$ and that $\tau$ is symmetric.

For $\tau\in \ccS_4^M$, the abelian variety $A_\tau$ has an automorphism $\phi$ 
of order three, in fact $\phi^2+\phi+1=0$,
defined by the following commutative diagram
{\renewcommand{\arraystretch}{1.3}
$$
\begin{array}{ccccccccc}
0&\longrightarrow&\ZZ^8&\stackrel{\Omega_\tau}{\longrightarrow}&\CC^4
&\longrightarrow &X_\tau&\longrightarrow&0\\
&&M\downarrow\phantom{M}&&
{}^tA^{-1}\downarrow\phantom{{}^tA^{-1}}&&
\phi\downarrow\phantom{\phi}&&\\
0&\longrightarrow&\ZZ^8&\stackrel{\Omega_\tau}{\longrightarrow}&\CC^4
&\longrightarrow &X_\tau&\longrightarrow&0
\end{array}
\qquad\Bigl(\mbox{use}\quad M\Omega_\tau\,=\,\Omega_\tau {}^t\!A^{-1}\,\Bigr)~.
$$
}% end stretch 

\subsection{Remark}\label{m22}
In the paper \cite{pmv} matrices $M_{p,q}$ were used
to investigate projective models of Picard modular varieties. The matrix $M_{2,2}$ is conjugated in $Sp(8,\ZZ)$ to $M$
(where each matrix has $2\times 2$ blocks):
$$
M_{2,2}\,:=\,\left(\begin{array}{cccc}
0& 0&-I&0\\
0&-I&0&I\\
I&0&-I&0\\
0&-I&0&0\end{array}\right)\,=\, TMT^{-1},\quad 
T\,:=\,
\left(\begin{array}{cccc}
I& I&0&-I\\
I&0&-I&I\\
0&I&I&-I\\
0&-I&0&I\end{array}\right)\quad(\in Sp(8,\ZZ))~.
$$
Thus $T$ maps $\ccS_4^M$ to $\ccS_4^{M_{2,2}}$, and this allows us to apply the results from \cite{pmv}.

\subsection{Abelian varieties of Weil type and the Hermite upper half space}
\label{weiltype}
Identifying $\CC^4$ with $T_0A_\tau$, the holomorphic tangent space to $A_\tau$ in 
$0\in A_\tau$, we see that the differential of $\phi$ is given by the matrix ${}^tA^{-1}$.
This matrix has $2$ eigenvalues $\omega$ and 2 eigenvalues $\overline{\omega}$ where
$\omega\in\CC$ is a primitive cube root of unity. Thus $A_\tau$ is an abelian variety of Weil type for the field $\QQ(\omega)$. The moduli space of such abelian varieties is isomorphic to $SU(2,2)/S(U(2)\times U(2))$ (cf.\ \cite[Prop. 5.5]{pmv}) and this symmetric domain is known as the Hermite upper half space $\ccH_{2,2}$.

\subsection{The theta functions}
We introduce the theta functions needed to find the projective model of the Picard modular variety.
For $\tau\in\ccS_4$ and $z\in\CC^4$ the classical theta function 
with characteristics
$[{}^\epsilon_{\epsilon'}]$, $\epsilon,\epsilon'\in\{0,1\}^4$, is defined by the series
$$
\theta[{}^\epsilon_{\epsilon'}](\tau,z)\,:=\,\sum_{m\in\ZZ^4}
e^{(m+\epsilon/2)\tau{}^t(m+\epsilon/2)\,+\,2(m+\epsilon/2){}^t(z+\epsilon'/2)}~.
$$\

The second order theta functions are the linear combinations of the $16$ functions $\theta[{}^\epsilon_0](2\tau,2z)$, $\epsilon\in\{0,1\}^4$. 
These functions define a holomorphic map
$X_\tau\rightarrow\PP^{15}$ which factors over the Kummer variety 
$X_\tau/\{\pm 1\}$.

We will study the map given by second order thetanulls:
$$
\Theta\,:\,\ccS_4\,\longrightarrow\,\PP^{15},\qquad 
\tau\,\longmapsto\,(\ldots: \theta[{}^\epsilon_0](2\tau,0):\ldots)~.
$$
This map factors over $\ccA_4(2,4):=\ccS_4/\Gamma_4(2,4)$, where $\Gamma_4(2,4)$ is a (normal) congruence subgroup of $Sp(8,\ZZ)$ (\cite[V.2, p.178]{Igusa}):
$$
\Gamma_g(2)\,:=\,\{M\in Sp(2g,\ZZ):\;M\equiv I\;\mbox{mod } 2\,\}~,
$$
$$
\Gamma_g(2,4)\,:=\,\left\{\left(\begin{array}{cc}a&b\\c&d\end{array}\right)
\in\Gamma_g(2):\,\mbox{diag}(a\,{}^t\!b)\equiv \mbox{diag}(c\,{}^t\!d)
\equiv 0\,\mbox{mod }2\,\right\}~.
$$
The image of $\Theta$ is a quasi projective variety.
The map $\Theta$ extends to a map from the Satake compactification 
$\overline{\ccA_4(2,4)}$ of $\ccA_4(2,4)$ to $\PP^{15}$
whose image is the closure of $\Theta(\ccS_4)$.

\subsection{The classical theta nulls and quadrics}\label{classqua}
We define quadrics in $\PP^{15}$ whose intersection with the Picard modular variety 
will provide information on special subloci.

The following formula relates the $\theta[{}^\epsilon_{\epsilon'}]$ 
to the second theta functions:
(cf.\ \cite[IV.1]{Igusa}, \cite[(3.3.2), (3.5.1)]{pmv}):
$$
\theta[{}^\epsilon_{\epsilon'}](\tau,z)^2\,=\,
\sum_\sigma (-1)^{(\sigma+\epsilon){}^t\epsilon'}\,
\theta[{}^\sigma_0](2\tau,0)\theta[{}^{\sigma+\epsilon}_{\;\;0}](2\tau,2z)~.
$$
Since $\theta[{}^\epsilon_{\epsilon'}](\tau,-z)=(-1)^{\epsilon{}^t\epsilon'}
\theta[{}^\epsilon_{\epsilon'}](\tau,z)$, 
the $136$ theta functions with $\epsilon{}^t\epsilon'=0$ are even functions of $z$, 
the remaining $120$ are odd. 

For an even characteristic $[{}^\epsilon_{\epsilon'}]$ and for $z=0$,
the formula above can be rewritten as
$$
\theta[{}^\epsilon_{\epsilon'}](\tau,0)^2\,=\,
Q[{}^\epsilon_{\epsilon'}](\ldots,\theta[{}^{\sigma}_{0}](2\tau,0),\ldots)~,
$$
where the functions $\theta[{}^\epsilon_{\epsilon'}](\tau,0)$ on $\ccS_4$ are called thetanulls, 
and where 
$$
Q[{}^\epsilon_{\epsilon'}]\,:=\,\sum_\sigma (-1)^{(\sigma+\epsilon){}^t\epsilon'}\,
X_\sigma X_{\sigma+\epsilon}~
$$ is a quadratic polynomial in the $16$ variables $X_\sigma$, with $\sigma\in(\ZZ/2\ZZ)^4$.
These quadratic polynomials define quadrics $Q[{}^\epsilon_{\epsilon'}]=0$ in $\PP^{15}$.

The points in the intersection of the closure of $\Theta(\ccS_4)$ and these quadrics 
correspond to (limits of) abelian varieties with vanishing thetanulls. 
Extending the results of  \cite[3.6,3.7, Lemma 3.8]{pmv} for $g=2,3$ to $g=4$,
we have the following table.
{\renewcommand{\arraystretch}{1.2}
$$
\begin{array}{|@{\hspace{10pt}}l|@{\hspace{10pt}}c|}
\hline
\mbox{moduli point}& \mbox{\# vanishing thetanulls}
\\ \hline
E\times A, \quad E\;\mbox{elliptic curve},\;A\; \mbox{abelian 3-fold} & 28\\ \hline
B_1\times B_2, \quad B_i\; \mbox{abelian surfaces} & 36\\ \hline
(\CC^\times)^2\times B,\quad\mbox{$B$ abelian surface}&96\\ \hline
(\CC^\times)^4, \quad \mbox{boundary point} & 120 \\ \hline
\end{array}
$$
}

\section{The Picard Modular variety $\ccX$}

\subsection{The action of $M$ on $\PP^{15}$}\label{actM}
The map $\Theta$ is equivariant for the action of $M$ on $\ccS_4$. In fact, from the definition of the theta functions it is obvious that
$$
\theta[{}^\epsilon_0](2(M\cdot\tau),0)\,=\,\theta[{}^{\epsilon A}_{\;0}](2\tau,0)~.
$$
For
$\rho,\sigma\in\ZZ^4$ one has
$\theta[{}^{2\rho+\sigma}_{\;\;\;0}](2\tau,0)=\theta[{}^\sigma_0](2\tau,0)$.
Thus the action of $M$ on these $16$ functions is determined by the action of $A$ on $(\ZZ/2\ZZ)^4$. 
We will write 
$$
\theta_i(\tau)\,:=\,\theta[{}^\epsilon_0](2\tau,0),\qquad\mbox{with}\quad
i\,:=\,\epsilon_12^3+\epsilon_22^2+\epsilon_32+\epsilon_4,\quad
\epsilon=(\epsilon_1,\ldots,\epsilon_4)\,\in\,\{0,1\}^4~.
$$
For example $\theta_{10}$ has characteristic $\epsilon=(1,0,1,0)$, which we simply write as $[1010]$. As $\epsilon A=(0,0,1,0)$ mod $2$ we see that $M$ maps $\theta_{10}$ to $\theta_2$.
With this convention, the action of $M$ is as follows:
$M\theta_0=\theta_0$, and $M$ has $5$ cycles of length three:

{\renewcommand{\arraystretch}{1.3}
$$
\begin{array}{rclclclclcrclclclcl}
\theta_1&\mapsto&\theta_4&\mapsto &\theta_5&\mapsto& \theta_1~,&\qquad
\theta_2&\mapsto&\theta_8&\mapsto &\theta_{10}&\mapsto &\theta_2~,\\
\theta_3&\mapsto&\theta_{12}&\mapsto &\theta_{15}&\mapsto& \theta_3~,&\qquad
\theta_6&\mapsto&\theta_{13}&\mapsto &\theta_{11}&\mapsto &\theta_6~,\\
\theta_7&\mapsto&\theta_{9}&\mapsto& \theta_{14}&\mapsto &\theta_7~.&&&&&&
\end{array}
$$
}% end stretch

Denote the coordinates on $\PP^{15}$ by $X_\epsilon$, $\epsilon\in(\ZZ/2\ZZ)^4$, or by $X_i$, with the relation between $\epsilon$ and $i$ as above. 
Then we define an action of $M$ on $\PP^{15}$ by a projective linear transformation as follows:
$$
M:\,\PP^{15}\,\longrightarrow\,\PP^{15},\qquad 
M(\ldots:X_\epsilon:\ldots)\,=\,(\ldots:X_{\epsilon A}:\ldots)~,
$$
so $M$ permutes the coordinates in the same way as it permutes the theta functions. 
Obviously,
we then have $\Theta(M\tau)=M\Theta(\tau)$, so the map $\Theta$ is equivariant for the actions of $M$ on $\ccS_4$ and $\PP^{15}$ respectively.

\subsection{The eigenspace $\PP^5$} \label{p5}
We are now interested in the image of $\ccS_4^M$ in $\PP^{15}$.
%it is a quasi projective variety which we denote by $X^0$.
As the map $\Theta$ is equivariant for $M$, 
this image must lie in an eigenspace for the action of $M$ on $\PP^{15}$.
In fact, as 
$\theta[{}^\epsilon_0](2(M\cdot\tau))=\theta[{}^{\epsilon A}_{\;0}](2\tau)$,
the image lies in the $M$-eigenspace, simply denoted by $\PP^5$, defined by the linear equations 
$$
\Theta:\ccS_4^M\,\longrightarrow\,\ccX\,\subset\,\PP^5\quad(\subset\PP^{15})~,
\qquad 
\PP^5:\quad X_\epsilon\,=\,X_{\epsilon A}\qquad (\epsilon\in(\ZZ/2\ZZ)^4)~. 
$$
This eigenspace $\PP^5\subset \PP^{15}$ can thus be parametrized by
$$
%(X_0:X_1:X_2:X_3:X_6:X_7)
X\,\longmapsto\,
(X_0:X_1:X_2:X_3:X_1:X_1:X_6:X_7:X_2:X_7:X_2:X_6:X_3:X_6:X_7:X_3)~,
$$
where $X=(X_0:X_1:X_2:X_3:X_6:X_7)$.

\subsection{The variety $\ccX$}\label{defx}
The fourfold $\Theta(\ccS_4^M)$ thus lies in this $\PP^5$, and we recall how one can find the equation for its closure $\ccX$. There are identities, valid for all $\tau\in\ccS_4$ 
and for some choice of signs, between even thetanulls of the form
$$
\prod_{i=0}^3\theta[{}^{\epsilon_{1i}00}_{\epsilon_{1i}'00}](\tau)\,\pm\,
\prod_{i=0}^3\theta[{}^{\epsilon_{2i}00}_{\epsilon_{2i}'00}](\tau)\,\pm\,
\prod_{i=0}^3\theta[{}^{\epsilon_{3i}00}_{\epsilon_{3i}'00}](\tau)\,\pm\,
\prod_{i=0}^3\theta[{}^{\epsilon_{4i}00}_{\epsilon_{4i}'00}](\tau)\,=\,0
$$
for suitable even $g=2$ characteristics $[{}^{\epsilon_{ni}}_{\epsilon_{ni}'}]$
(cf.\ \cite[4.4]{pmv},\cite[1.8]{order5}). Actually the characteristics given in \cite[1.8]{order5} contain a misprint (one of them is odd!), below are two sets of four $g=2$ characteristics for which we have identities as above:
$$
[{}^{00}_{00}],\;[{}^{00}_{10}],\;[{}^{10}_{00}],\;[{}^{11}_{11}],\qquad
\mbox{and}\quad
[{}^{01}_{00}],\;[{}^{01}_{10}],\;[{}^{10}_{00}],\;[{}^{10}_{01}]~.
$$
By taking the product of the eight expression on the left hand side for all choices of signs, one obtains a polynomial in the $\theta[{}^\epsilon_{\epsilon'}]^2$. Using the formula from section  \ref{classqua}, this can be written as a polynomial, of degree $32$, in the sixteen second order theta constants $\theta[{}^\sigma_0](2\tau,0)$.  The zero locus of this polynomial in $\PP^{15}$ contains the image of $\ccS_4$, and thus restricting it to the $\PP^5$ gives an equation for $\ccX$.  
Taking the GCD of two such equations of degree $32$, 
we found that the image is defined by a polynomial $F$ of degree $10$.
Thus 
$$
\ccX\,=\,\overline{\Theta(\ccS_4^M)}\,=\,
\Theta\big(\overline{\ccA_4(2,4)}\big)\,\cap\,\PP^5\,=\,
Z(F)
$$
where $\overline{\ccA_4(2,4)}$ is the Satake compactification of $\ccA_4(2,4)=\ccS_4/\Gamma_4(2,4)$ and
$Z(F)$ is the zero locus of $F$ in $\PP^5$.

The polynomial $F$ defining $\ccX$ is homogeneous of degree $10$ in the six variables $X_0$, $X_1$, $X_2$, $X_3$, $X_6$, $X_7$, which give the coordinate functions on $\PP^5$, and it has $147$ terms. 
%This polynomial 
It is symmetric in $X_1,\ldots,X_7$, and can be written as
$$
F\,:=\,F_{10}\,-\,X_0X_1X_2X_3X_6X_7F_4
$$
where the homogeneous polynomials $F_4,F_{10}$ of degree $4$ and $10$ respectively,
are given by
{\renewcommand{\arraystretch}{1.3}
$$
\begin{array}{rcl}
F_4&:=&-6 S_1^2 + 16 S_2 + 4 S_1 X_0^2 + 2 X_0^4~, \\
F_{10}&:=&S_1 S_2^2-3 S_1^2 S_3+12 S_1 S_4-48 S_5+(-S_2^2+2 S_1 S_3+4 S_4) X_0^2+S_3 X_0^4~,
\end{array}
$$
}
where $S_i(X_1,\ldots,X_7):=s_i(X_1^2,\ldots,X_7^2)$ and $s_i$ is the degree $i$
elementary symmetrical function in the five variables $X_1,\ldots,X_7$.

\subsection{The singular locus of $\ccX$}\label{singX}
The polynomial $F$ defining $\ccX$ is rather complicated. 
We relied on Magma to show that the singular locus $Sing(\ccX)$ of $\ccX$ 
has dimension two and degree $320$. 
It was then not hard to find $120$ quadratic surfaces and $80$ planes in $Sing(\ccX)$, so we accounted for all two-dimensional components of $Sing(\ccX)$. 
The general points of these components all correspond to decomposable ppav's, see Propositions \ref{psing22} and \ref{psing31} below. 
These points are in fact quotient singularities in $\ccA_4(2,4)$
%$\ccS_4/\Gamma_4(2,4)$. 
We do not know if there are components of lower dimension in $Sing(\ccX)$.

\

\section{Automorphisms of $\ccX$}

\subsection{The Weyl group $W(E_6)$ and $Aut(\ccX)$}
The subgroup of $Sp(8,\ZZ)$ of elements which map $\ccS_4^M$ into itself 
acts by projective transformations on the eigenspace $\PP^5$ 
of $M$ in $\PP^{15}$ and maps $\ccX=\overline{\Theta(\ccS_4^M})$ into itself.
We recall the results from \cite{pmv} on this subgroup and we show that 
it acts as the Weyl group of the root system $W(E_6)$ on $\PP^5$.

\subsection{Centralizers and normalizers}\label{cennor}
The normalizer of the subgroup 
$\langle M\rangle=\{I,M,M^{-1}\}$ in $Sp(8,\ZZ)$
is the subgroup 
$$
N_M\,:=\,\{A\,\in\,Sp(8,\ZZ):\,AMA^{-1}\,=\,M^{\pm 1}\,\}~,
$$
whereas the centralizer of $\langle M\rangle$ is the subgroup:
$$
C_M\,:=\,\{A\,\in\,Sp(8,\ZZ):\,AM\,=\,MA\,\}~.
$$
An element $A\in N_M$ will either permute $M,M^{-1}$ or it will fix 
both of them and in that case $M\in C_M$. 
Thus the index of $C_M$ in $N_M$ is either one or two.

The index $[N_M:C_M]=2$ since $M_BMM_B^{-1}=M^{-1}$, where
$$
B\,:=\,\left(\begin{array}{cc}0& I\\I&0 \end{array}\right)\quad(\in GL(4,\ZZ)),
\qquad
M_B\,:=\,\,\left(\begin{array}{cccc}B&0\\0&B \end{array}\right)
\quad(\in Sp(8,\ZZ))~,
$$
here $I$ is the $2\times 2$ identity matrix. 

\subsection{The action of the normalizer}
Notice that the normalizer $N_M$ maps the Hermite upper half space $\ccS_4^M$ into itself:
if $\tau\in \ccS_4^M$ then
$M\tau=\tau$, hence also $M^{-1}\tau=\tau$ and: 
$$
M(A\tau)\,=\,A(M^{\pm 1}\tau)\,=\,A\tau\qquad(A\,\in\, N_M,\;\tau\,\in \, \ccS_4^M)~.
$$ 
So we have biholomorphic maps
$$
A\,:\,\ccS_4^M\,\longrightarrow\,\ccS_4^M\qquad (\,A\,\in\,N_M\,)~.
$$

In the projective representation of $Sp(8,\ZZ)$ on $\PP^{15}$, any $A\in N_M$
will permute the eigenspaces for the action of $M$ on $\PP^{15}$. 
As $M$ has a unique $5$-dimensional eigenspace $\PP^5$, 
we must have $A(\PP^5)\subset\PP^5$ and thus
$$
A\,:\,\PP^5\,\longrightarrow\,\PP^5\qquad (\,A\,\in\,N_M\,)~.
$$
Obviously, $\ccX$ is mapped into itself and thus we have a homomorphism
$N_M\rightarrow Aut(\ccX)$.

\subsection{A Hermitian form}\label{HM}
The action of $M\in Sp(8,\ZZ)$, with $M^2+M+I=0$, on $\ZZ^8$ defines
the structure of $\ZZ[\omega]\cong\ZZ[x]/(x^2+x+1)$ module on $\ZZ^8$, 
where $\omega$ acts as $M$.
Using the alternating form on $\ZZ^8$ defined by the matrix $E$, 
which we denote also by $E$,
we define a Hermitian form on $\ZZ^8$
with values in $\ZZ[\ww]$ as follows
(cf. \cite[Lemma 6.2]{pmv} and the proof of Proposition \ref{propher} below):
$$
 H_M(x,y)\, := \,E(x,My) -\omega E(x,y)~.
$$

Let  $e_i$ be the $i$-th standard basis vector of $\ZZ^8$.
Then we have the Gram matrix: 
$$
\big(H_M(f_i,f_j)\big)=\mbox{diag}(1,1,-1,-1),\qquad
f_i:\,=\,e_i+e_{i+4},\qquad
(i,j\,=\,1,\ldots,4)~,
$$
so the $f_i$ are an orthogonal $\ZZ[\ww]$-basis of $\ZZ^8$ and
the signature of $H_M$ is $(2,2)$.

The quotient ring $\ZZ[\omega]/2\ZZ[\omega]\cong (\ZZ/2\ZZ)[x]/(x^2+x+1)$ 
is isomorphic to the finite field $\FF_4$. The Hermitian form $H_M$ on $\ZZ^8$
defines a Hermitian form on $\FF_4^4\cong (\ZZ[\omega]/2\ZZ[\omega])^{4}$
by reduction modulo $2$.

We recall the following results.

\subsection{Proposition}\label{propher}
\begin{enumerate}
\item
The centralizer $C_M(\RR)$ of $M$ in $Sp(8,\RR)$ is isomorphic to $U(H_M)\cong U(2,2)$, the unitary group of the 
($\RR$-bilinear extension of the) Hermitian form $H_M$
on $\ZZ^8\otimes_\ZZ\RR$.
\item 
The reduction modulo $2$ map induces a surjective homomorphism
from $C_M\subset Sp(8,\ZZ)$ onto $U(4,\FF_4)$. 
\item
The center of $U(4,\FF_4)$ is cyclic of order $3$ and is generated by the scalar multiplication $v\mapsto \ww v$. The quotient group $PU(4,\FF_4)$ is a finite simple group of order $25920$. 
\item 
The reduction modulo 2 map followed by the quotient by $<M>$-map is a surjective
homomorphism of $N_M$ onto $W(E_6)$. 
In particular,  $PU(4,\FF_4)$ is isomorphic to a subgroup of index two of the 
Weyl group $W(E_6)$ of the root system $E_6$.
\item The map $\Theta:\ccS_4^M\rightarrow\PP^5$
factors over $\ccS_4^M/C_M(2)$ where 
$$
C_M(2)\,:=\,\{A\,\in\,C_M\,:\;A\,\equiv\,I\,\mbox{mod}\;2\,\}\;=\;C_M\,\cap\,\Gamma_g(2)~.
$$
\end{enumerate}

\ts
The isomorphism $C_M(\RR)\cong U(H_M)$ is proven in \cite[Proposition 5.5.2]{pmv}.
The reduction map is studied in \cite[Lemma 6.2.3]{pmv}, 
the structure of $U(4,\FF_4)$ and its relation with $W(E_6)$
can be found in the Atlas \cite{Atlas}.
The last item is obtained from 
\cite[Proposition 6.4]{pmv}.
For completeness sake we check that $H_M$ is Hermitian.
Using that $E$ is alternating, $E(Mx,My)=E(x,y)$, $M^2=-I-M$, 
and $1+\omega=-\omega^2=:\overline{\omega}$ we have
$$
\begin{array}{rcl}
H_M(y,x) &=& E(y,Mx)\, -\,\omega E(y,x) \\
&=& -E(Mx,y) + \omega E(x,y) \\
&=& -E(M^2x,My) +\omega E(x,y) \\
&=& E(x,My)\,+\,E(Mx,My)\,+\,  \omega E(x,y)\\
&=& E(x,My)\,-\,\omega^2 E(x,y)\\
&=&\overline{H_M(x,y)}^{\phantom{x}}~.
\end{array}
$$
We also have 
$$
\begin{array}{rcl}
H_M(Mx,y)&=&E(Mx,My)\,-\,\omega E(Mx,y)\\
&=&E(x,y)\,-\,\omega E(M^2x,My)\\
&=&E(x,y)\,+\,\omega(E(x,My)+E(Mx,My))\\
&=&-\omega^2E(x,y)\,+\,\omega E(x,My)\\
&=&\omega H_M(x,y)~,
\end{array}
$$
so $H_M$ is indeed $\ZZ[\omega]$-linear in the first variable.
\qed

\

\subsection{The invariant quadric}\label{invqua}
The action of $N_M$ on $\PP^5$ thus induces an action of $W(E_6)$ on $\PP^5$. 
We will show in Proposition \ref{actwe6}
that this action is obtained from the standard reflection representation on the root lattice $R(E_6)$ by complexifying and projectivization:
$$
\PP^5\,\cong\,\PP(\,R(E_6)\,\otimes_\ZZ\,\CC\,)~.
$$

Among the $136$ quadratic forms $Q_m$, 
only one is $M$-invariant. It is
$Q[{}^0_0]=\sum X_\sigma^2$,
which restricts to (cf.\ Section \ref{p5})
$
X_0^2+3(X_1^2+X_2^2+X_3^2+X_6^2+X_7^2)
$ 
on $\PP^5$.
The bilinear form defined by this quadratic form, 
for convenience multiplied by a scalar,  will be denoted by $b$:
$$
b(X,Y)\,:=\,\mbox{$\frac{1}{3}$}X_0Y_0\,+\,X_1Y_1\,+\,X_2Y_2\,+\,X_3Y_3\,+X_6Y_6\,+\,X_7Y_7~.
$$
We will use  $b$ to define the inner product on the root lattice $E_6$. 

\subsection{The root system $E_6$}\label{E6}
The root system $E_6$ is defined by the Dynkin diagram:
$$
\put(-120,0){
\DynkinEEEE{\alpha_1}{\alpha_2}{\alpha_3}{\alpha_4}{\alpha_5}{\alpha_6}
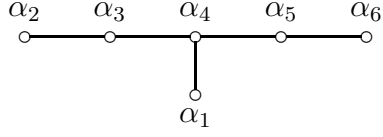
{\mbox{The Dynkin diagram of $E_6$}
}}
$$
(so $b(\alpha_i,\alpha_j)=0,-1,2$ if $\alpha_i$ and $\alpha_j$ are not connected,  are connected or $i=j$ respectively). 
The following basis of the simple roots of the root system $E_6$ will be related to action of $N_M$ on $\PP^5$:
$$
\begin{array}{rclrclrcl}
\alpha_1&:=&( 0,-1,-1,0,0,0)&\quad 
\alpha_3&:=&(0,1,-1,0,0,0),&\quad 
\alpha_5&:=&(0,0,0,1,-1,0),\\
\alpha_2&:=&(3,-1,1,1,1,1)/2,&\quad 
\alpha_4&:=&(0,0,1,-1,0,0),&\quad 
\alpha_6&:=&(0,0,0,0,1,-1)~.
\end{array}
$$
The root lattice of $E_6$ is $R(E_6)=\oplus_{i=1}^6\ZZ\alpha_i$.
A root $\alpha$ of $E_6$ defines a reflection on $\CC^6:=R(E_6)\otimes_\ZZ\CC$:
$$
s_\alpha:\,\CC^6\,\longrightarrow\,\CC^6,\qquad
s_\alpha(X)\,=\,X\,-\,\frac{2b(X,\alpha)}{b(\alpha,\alpha)}\alpha\;=\;
X\,-\,b(X,\alpha)\alpha
$$
where we used that $b(\alpha,\alpha)=2$. The reflections in the simple 
roots generate the Weyl group $W(E_6)$ of $E_6$, which is a finite group of order
$51840$.

The following proposition gives explicit matrices in $N_M\subset Sp(8,\ZZ)$ whose action on $\PP^5$ generates the group $W(E_6)$.

\subsection{Proposition}\label{actwe6}\label{actionMB}
Let $M_B\in N_M$  be as in Section \ref{cennor} and define $M_d,M_e,M_f\in N_M$ by:
% MB1, Mg4
$$
M_d\,:=\,\left(\begin{array}{rrrrrrrr}
 1& 0& 0& 0& 2& 0& -1& 0\\ 
 0& 1& 0& 0& 0& 0& 0& 0\\ 
 0& 0& 1& 0& -1& 0& 2& 0\\ 
 0& 0& 0& 1& 0& 0& 0& 0\\
 0& 0& 0& 0& 1& 0& 0& 0\\ 
 0& 0& 0& 0& 0& 1& 0& 0\\ 
 0& 0& 0& 0&0& 0& 1& 0 \\ 
 0& 0& 0& 0& 0& 0& 0& 1
 \end{array}\right)~,\quad
M_e\,:=\,\left(\begin{array}{rrrrrrrr}
-1&0&0&0&0&1&0&-2\\
0&1&0&0&-1&0&-1&0\\
0&0&-1&0&0&1&0&1\\
0&0&0&1&2&0&-1&0\\
0&0&0&0&-1&0&0&0\\
0&0&0&0&0&1&0&0\\
0&0&0&0&0&0&-1&0\\
0&0&0&0&0&0&0&1
\end{array}\right)
$$

$$
M_f\,:=\,\left(\begin{array}{cccc}0&B_f\\-B_f&0 \end{array}\right)
\quad\mbox{with}\quad 
B_f\,:=\,\left(\begin{array}{cc}I& 0\\0&-I \end{array}\right)
$$
where $I$ is the $2\times 2$ identity matrix.

The action of these elements in $N_M$ on $\PP^5$ is induced by the 
following linear transformations in $W(E_6)$:
{\renewcommand{\arraystretch}{1.3}
$$
\begin{array}{cccc}
M_B\,:&(X_0,X_1,X_2,X_3,X_6,X_7)&\longmapsto&(X_0,X_1,X_2,X_3,X_7,X_6)~\\
M_d\,:&(X_0,X_1,X_2,X_3,X_6,X_7)&\longmapsto&(X_0,X_1,-X_2,-X_3,-X_6,-X_7)~,\\
M_e\,:&(X_0,X_1,X_2,X_3,X_6,X_7)&\longmapsto&(X_0,X_1,X_2,-X_3,-X_6,X_7)~,
\end{array}
$$
} %end stretch
$$
M_f\,:\,X\, \longmapsto\,\mbox{$\frac{1}{4}$}\left(\begin{array}{rrrrrr}
1&3&3&3&3&3\\
1&-1&3&-1&-1&-1\\
1&3&-1&-1&-1&-1\\
1&-1&-1&3&-1&-1\\
1&-1&-1&-1&-1&3\\
1&-1&-1&-1&3&-1
\end{array}\right)X~,
$$
where $X=(X_0,\ldots,X_7)$ is viewed as a column vector.
In particular, $M_B$ acts as $s_{\alpha_6}$. Moreover, $W(E_6)$ is generated by these four linear transformations on $\CC^6$.

\ts
The action of $M_B$ on $\PP^{15}$ is very similar to the one of $M$ which we found in 
Section \ref{actM}, it is simply $M_B:X_\sigma\mapsto X_{\sigma B}$.
Thus under the action $M_B$ the coordinates $X_0,X_5,X_{10},X_{15}$ are fixed
and the remaining ten are permuted as follows:
$$
X_1\,\leftrightarrow\,X_4,\quad
X_2\,\leftrightarrow\,X_8,\quad
X_3\,\leftrightarrow\,X_{12},\quad
X_6\,\leftrightarrow\,X_{9},\quad
X_7\,\leftrightarrow\,X_{13},\quad
X_{11}\,\leftrightarrow\,X_{14}~.
$$
On the eigenspace $\PP^5$ of $M$, parametrized as in Section \ref{p5},
we then have
$$
M_B:\,\PP^5\longrightarrow\,\PP^5,\qquad
M_B:\,(X_0:X_1:X_2:X_3:X_6:X_7)\,\longmapsto\,(X_0:X_1:X_2:X_3:X_7:X_6)~,
$$
for example $M_B:X_6\mapsto X_9$, but on $\PP^5$ we have $X_7=X_9=X_{14}$.

The action of $M_d$ and $M_f$ is easy to compute using the series expansion of the theta functions. To find the action of $M_f$, one can use that 
$M_f=\mbox{diag}(I,-I,I,-I)E$
where we take $2\times 2$ diagonal blocks and $E$ is as in Section \ref{auto}.
The action of $E$ on the second order thetanulls is well-known, on $\PP^{15}$ it is given by $X_\sigma\mapsto\sum_\rho(-1)^{\sigma{}^t\!\rho}X_\rho$. 
By restricting to $\PP^5$ we find the matrix in $W(E_6)$.
A Magma computation shows that the matrices which define the action on $\PP^5$ indeed generate $W(E_6)$.
\qed

\subsection{Invariants of $W(E_6)$}
The equation $F$ for $\ccX$, which we determined in Section \ref{defx}, is an invariant for the $W(E_6)$-action. The ring of $W(E_6)$-invariant polynomials in $\CC[X_0,\ldots,X_7]$
is generated by invariants $I_k$ of degree $k$, for $k=2,5,6,8,9,12$.
These  are defined as the sum of the $k$-th powers of the
hyperplanes perpendicular to the 27 vectors
$v_1,\ldots,v_{27}$ in the $W(E_6)$-orbit of the vector
$v_1=(1,0,\ldots,0)$, so $ I_k = \sum_{i=1}^{27} b(X,v_i)^k$. 
For example, $I_2=(3/2)b(X,X)$.
A computation shows that, with $c=-2/675$, the polynomial $F$ defining $\ccX$ can be written
as
$$
F\,=\, c(11520 I_8 I_2\, -\,4160 I_6 I_2^2\, -\,4608 I_5^2\, + \,25 I_2^5)~.
$$

\subsection{The quotient $\overline{\ccX}:=\ccX/W(E_6)$}
The quotient variety $\PP^5/W(E_6)$ is the weighted projective space 
$W\PP^5:=W\PP^5(2,5,6,8,9,12)$ 
and the quotient map is given by the invariants $I_k$. 
The projection of $\ccX/W(E_6)\subset W\PP^5$ onto $W\PP^4:=W\PP^4(2,6,8,9,12)$ 
is then a 2:1 branched cover with covering involution given by 
$I_5\mapsto -I_5$.
Note that the branch locus is reducible, 
one component is defined by $I_2=0$, the other by
$11520 I_8  -4160 I_6 I_2 + \,25 I_2^4=0$.

\

\section{Decomposable abelian varieties}\label{subX}

\subsection{The two cases}\label{redprop}
We consider the case that the abelian fourfold $A_\tau$, for $\tau\in \ccS_4^M$, is a product of lower dimensional ppav's. We will show the following:
\begin{enumerate}
\item The closure of the locus in $\ccX$ of ppav's which are products of two abelian surfaces consists of 120 (smooth) quadric surfaces
(see Section \ref{sing22}).

\item The closure of the locus in $\ccX$ of ppav's which are products of an elliptic curve and an abelian threefold consists of 80 projective planes, we refer to these as Hesse planes
(see Section \ref{sing31}).
\end{enumerate}
In both cases, the surfaces which parametrize these products lie in $Sing(\ccX)$.

\subsection{The products of abelian surfaces}\label{sing22}
The moduli space of two dimensional ppav's 
with an automorphism of order three (of type $(1,1)$) 
and a level two structure has a model which is a $\PP^1\subset \overline{\ccA_2(2,4)}$ (\cite[Theorem 8.4]{pmv}).
There are 3 points on this $\PP^1$ where the abelian surface degenerates 
to $(\CC^\times)^2$ and there are two points where it decomposes into $E_3^2$,
where the elliptic curve $E_3$ is defined as
$$
E_3\,:=\,\CC/(\ZZ+\omega\ZZ)~.
$$
The product of two such abelian surfaces is a fourfold of the type we consider here, 
so we expect to see copies of $\PP^1\times\PP^1$ inside $\ccX$. In fact, the quadric $Q_{22}$ in the next proposition parametrizes such products.

The roots $\alpha_5=(0,0,0,1,-1,0)$, $\alpha_6=(0,0,0,0,1,-1)\in E_6$ 
define the hyperplanes $\alpha_5^\perp:\,X_3=X_6$ and 
$\alpha_6^\perp:\,X_6=X_7$.
Thus the projective 3-space $Z$ in the following proposition is $Z=\alpha_5^\perp\cap\alpha_6^\perp$. These two roots
span the root system$\{\pm\alpha_5,\pm\alpha_6,\pm(\alpha_5+\alpha_6)\}$ in $E_6$
which is isomorphic to $A_2$. 
There are $120$ such subsystems and $W(E_6)$ acts transitively on them, so we get $120$ quadrics like $Q_{22}$ in $\ccX$.

\subsection{Proposition}\label{psing22}
Let $Z\subset\PP^5$ be the projective 3-space defined by
$$
Z\,:\quad X_3\,=\,X_6\,=\,X_7\quad (\subset\,\PP^5)~.
$$
The intersection of $Z$ with $\ccX\subset\PP^5$ has two irreducible components,
$\ccX\,\cap\,Z\;=\;Q_{22}\,\cup\,S_{22}$,
where $Q_{22}$ is a quadric which lies in $Sing(\ccX)$ and $S_{22}$ is a sextic surface.

The quadric $Q_{22}$ parametrizes products of abelian surfaces, each with an automorphism of order three of type $(1,1)$.
The surface $S_{22}$, which is birationally isomorphic to a K3 surface, 
parametrizes abelian fourfolds which are isogeneous to a product of abelian surfaces. 

\ts
An explicit computation shows that the restriction of $F$ to $Z$ factors as
$$
F(X_0,X_1,X_2,X_3,X_3,X_3)\,=\,q_{22}^2f_{22},\qquad
q_{22}\,:=\,X_0X_3\,-\,X_1X_2~.
$$
As $Q_{22}$ is a smooth quadric, it is isomorphic to $\PP^1\times\PP^1$.
A parametrization is given by
$$
\PP^1\times\PP^1\,\longrightarrow\,Q_{22}\;(\subset Z),\quad
\Big((s:t),(u:v)\Big)\,\longmapsto\,(X_0:X_1:X_2:X_3)\,=\,(su:sv:tu:tv)~.
$$
It is now straightforward to check that all partial derivatives of $F$ vanish on $Q_{22}$,
hence $Q_{22}\subset \mbox{Sing}(\ccX)$.
We checked that exactly $36$ of the quadrics $Q_m$ vanish in a general point of $Q_{22}$, so such a point corresponds to a product of two abelian surfaces.

To find a subdomain in $\ccS^M_4$ which maps to $Q_{22}$, 
we consider the following matrix:
$$
M_{pr}\,:=\,\left(\begin{array}{cc}A_p&0\\0&{}^tA_p^{-1} \end{array}\right)
\quad(\in Sp(8,\ZZ))~,\qquad
A_p\,:=\,\left(\begin{array}{cccc}
-1&0&-1&0\\
0&1&0&0\\
1&0&0&0\\
0&0&0&1 \end{array}\right)~,
$$
so $A_p$ acts as $A$ (the $(1,1)$-block of $M$) on the first and the third coordinate, but it acts as the identity on the second and fourth coordinate. One has $M_{pr}M=MM_{pr}$, so $M_{pr}\in C_M$ acts on the Hermite upper half space $\ccS_4^M$. 
The fixed point set of $M_{pr}$ in $\ccS_4^M$ consists of the
block matrices (with indices $1,3$ and $2,4$ respectively) 
of the period matrices of abelian surfaces with an automorphism of order three.
In particular, $\ccS_4^{M,M_{pr}}$ parametrizes products of ppav's.
The map $\Theta$ will map this 2-dimensional domain 
to an eigenspace of the action of $M_{pr}$ on $\PP^5$.
The  action of $M_{pr}$ on $\PP^5$ is given by $X_\sigma\mapsto X_{\sigma A_p}$,
and one easily checks that $M_{pr}$ fixes $Z$ pointwise and that it has two other isolated fixed points in $\PP^5$.
For dimension reasons, we get that
$\Theta(\ccS_4^{M,M_{pr}})\subset Z$, 
and the closure of the image must be $Q_{22}$.

Recall that the reflection defined by $\alpha_6$ is induced by $M_B$, 
as in Section \ref{actionMB}.
The  matrices $M_B$ and $M_{pr}$ generate a subgroup of order six in $Sp(8,\ZZ)$ 
which is isomorphic to $W(A_2)\cong S_3$. 
One reason for this is that the an abelian surface with an automorphism 
of order three of type $(1,1)$ actually admits $S_3$ as automorphism group 
(the case III in \cite[Section 11.7]{BL}).

Now we discuss the other component $S_{22}$. To see a subdomain of the 
Hermite upper half space which maps to $S_{22}$ we consider the matrix
$$
M_{ip}\,:=\,\left(\begin{array}{cc}A_p&B_{ip}\\0&{}^tA_p^{-1} \end{array}\right)
\quad(\in Sp(8,\ZZ))~,\qquad
B_{ip}\,:=\,\left(\begin{array}{cccc}
0&-2&0&0\\
2&0&0&0\\
0&2&0&-2\\
-2&0&2&0 \end{array}\right)~,
$$
where $A_p$ is as above. In particular, $M_{pr}$ and $M_{ip}$ have the same image in $\Gamma_4/\Gamma_4(2,4)$ and thus $Z$ is also an eigenspace of the 
action of $M_{ip}$ on $\PP^5$.
The matrices $M_{ip}$ and $M_B$
also generate a subgroup of $Sp(8,\ZZ)$ isomorphic to $S_3$.
However,  the fixed point loci $\ccS_4^{M,M_{pr}}$ and $\ccS_4^{M,M_{ip}}$
are not conjugate under the action of $Sp(8,\ZZ)$. In fact, consider
the sublattices 
$\Lambda_{pr}:=\ker(M_{pr}-I)$ and $\Lambda_{ip}:=\ker(M_{ip}-I)$,
both isomorphic to $\ZZ^4$. 
The alternating form $E$ restricts to an alternating form with determinant $1$ on $\Lambda_{pr}$, 
but its restriction to $\Lambda_{ip}$ has determinant $9$, which implies the matrices cannot be conjugate in $Sp(8,\ZZ)$. 
As $\Theta$ also maps $\ccS_4^{M,M_{ip}}$ to $Z$, 
the image $\Theta(\ccS_4^{M,M_{ip}})$ 
must be the other component $S_{22}$ of the intersection of $\ccX$ with $Z$.

The surfaces $Q_{22}$ and $S_{22}$ intersect along a curve of degree $12$, which is the union of six lines, each with multiplicity two. 
These lines are also the singular locus of $S_{22}$. 
A better model of $S_{22}$ can be obtained as the image of the birational map
$$
S_{22}\,\longrightarrow\,\PP^4,\qquad
(X_0:\ldots:X_3)\,\longmapsto (X_0q_{22}(X):\ldots:X_3q_{22}(X):r_3(X))~,
$$
with
$$
r_3(X)\,:=\,X_3(X_0-X_1-X_2+X_3)(X_0+X_1+X_2+X_3)~.
$$
All coordinate functions vanish on the six lines and are homogeneous of degree $3$.
The birational inverse of this map is induced by the projection on the first four coordinates.
The image of this map was found with Magma. 
It is a complete intersection of a quadratic and a cubic hypersurface. 
The image is smooth except for $9$ ordinary double points.
Thus the minimal model of the image, and hence of $S_{22}$, is a K3 surface. 
\qed

\subsection{The Hesse planes} \label{sing31}
In \cite[ Theorem 8.5]{pmv}, it was shown that the moduli space of three dimensional ppav's with an automorphism of order three (of type $(2,1)$)
and a level two structure has a projective model which is a $\PP^2$.  
Taking the product of such a threefold with
the elliptic curve $E_3$ (see Section \ref{sing22}),
we obtain an abelian fourfold of the type we consider here.
The decomposable ppavs, products of a abelian surface with 
the elliptic curve $E_3$, form a configuration of 12 lines. 
These are the 12 lines in the four reducible curves in the Hesse pencil $x^3+y^3+z^3+\lambda xyz=0$ (cf.\ \cite{AD}). 
A Hesse plane will be a copy of a $\PP^2$ with a Hesse pencil.
Thus we expect to find such Hesse planes inside $\ccX$.

\subsection{Proposition}\label{psing31}
There is a unique conjugacy class $C$ in $W(E_6)$
consisting of 80 elements, each of order three, 
whose characteristic polynomial in the six dimensional representation is $(x^2+x+1)^3$. 

Each $g\in C$ has two eigenspaces in $\PP^5$, both are planes 
$\PP^2\subset \ccX$. Moreover, each of the two planes lies in $Sing(\ccX)$
and in this way we get $80$ planes in $Sing(\ccX)$.
Each plane parametrizes products of an abelian threefold and an elliptic curve.

\ts
The conjugacy class $C$ can be found from 
\cite[Table II, p.104]{Frame} or \cite{Atlas}.
The $80$ elements in $C$ come in pairs, $g$, $g^2$, which have the same eigenspaces.
The group $W(E_6)$ thus acts transitively on the set of $2\cdot 40=80$
eigenspaces of the elements of $C$.

To be explicit, here is one element $g_3\in C$ 
and one of its eigenspaces $W_3$:
$$
g_3\,:=\,\frac{1}{2}\left(\begin{array}{rrrrrr}
-1&0&0&0&0&-3\\
0&-1&-1&1&1&0\\
0&1&-1&1&-1&0\\
0&-1&-1&-1&-1&0\\
0&-1&1&1&-1&0\\
1&0&0&0&0&-1
\end{array}\right),\qquad
W_3\,:=\,\left(\begin{array}{rrr}
3&0&0\\
0&1&0\\
0&0&1\\
0&\ww&\;-\ww^2\\
0&\;-\ww^2&-\ww\\
-1-2\ww&0&0
\end{array}\right)~,
$$
where the columns of $W_3$ span the eigenspace of $g_3$ with eigenvalue $\ww$. 
We verified that $W_3\subset Sing(\ccX)$ and that $28$ of the quadrics $Q_m$ vanish on $W_3$. 
Hence $W_3$ parametrizes products of an abelian threefold and an elliptic curve. 
The other quadrics $Q_m$ intersect each plane in two lines.
\qed

\subsection{Remark}
The centralizer of $g_3$ in $W(E_6)$ acts on the eigenspace $\PP W_3$
as in the proof of the proposition. 
With Magma we found that it is a group of order 
$648$ which coincides with the group denoted by $\overline{G}_{216}$ in \cite{AD} acting on the Hesse pencil.
The $W(E_6)$-invariants $I_k$, for $k=2,5,8$ restrict to zero on $\PP W_3$, whereas
$I_6,I_9,I_{12}$ restrict to invariants of $\overline{G}_{216}$.
We checked that the restrictions of $I_6^2$ and $I_{12}$ are linearly independent degree $12$ polynomials.
The ring of invariants of $\overline{G}_{216}$ is thus generated by
the restrictions of $I_6,I_9,I_{12}$ (cf.\ \cite{AD}).

\

\section{Some fixed point loci in $\ccX$}

\subsection{Fixed points}
In this section, we consider the fixed point loci in $\ccX$ of one reflection
and of two commuting reflections in $W(E_6)$. The case of two non-commuting reflections was already described in Proposition \ref{psing22}.
We show that the Hessian ${\mathcal W}_{10}$ of Igusa's quartic threefold
is an arithmetic quotient. This was conjectured by Hunt (\cite[p.7-8]{Hunt}),
based on an analogy with the Nieto quintic.
We also find a projective model of a two-dimensional moduli space 
of abelian fourfolds with an automorphism of order $12$.

\subsection{Proposition} \label{fixMB}
Let $H_\alpha\subset\PP^5$ be the projectivization of the reflection hyperplane defined by a root $\alpha\in E_6$. Then the intersection
of $H_\alpha$ with $\ccX$ is an irreducible 3-fold of degree $10$,
which is ${\mathcal W}_{10}$, the Hessian of the Igusa quartic:
$$
{\mathcal W}_{10}\,:=\,H_\alpha\cap \ccX~.
$$
This 3-fold parametrizes abelian fourfolds of Weil type which are isogeneous to the selfproduct of an abelian surface.

\ts
As $W(E_6)$ acts transitively on the roots of $E_6$, it suffices to consider
the case that $\alpha:=(-3,1,1,1,1,1)/2$.
Then $H_\alpha\subset \PP^5$ is defined by the linear equation 
$$
H_\alpha\,:=\,\{\ccX\,\in\,\PP^5\,:\; b(\alpha,X)\,=\,0\;\}\,=\,
\{X\,\in\,\PP^5\,:\;-X_0\,+\,X_1\,+\,\ldots\,+\,X_7\,=\,0\;\}~.
$$

The equation of the intersection 
$H_\alpha\cap \ccX$ 
is thus 
$F(X_1+\ldots+X_7,X_1,\ldots,X_7)=0$. 
This polynomial is quite complicated, it has $591$ terms. 
However, an explicit computation shows that it can also be obtained as follows.
Let ${\mathcal I}_4$ be Igusa's quartic threefold in $\PP^4$, with coordinates $X_1,X_2,X_3,X_6,X_7$, which is defined by the equation 
$$
{\mathcal I}_4:\quad G\,:=\,s_2^2-4s_4\,=\,0,\qquad(\subset\PP^4)
$$
where the $s_i$ are the elementary symmetric functions of degree $i$ in these variables (cf. \cite[Theorem 5.2, Theorem 4.1]{vdG}, \cite[Section 3.3]{Hunt}).
Then we have, for a non-zero constant $c$:
$$
F(X_1+\ldots+X_7,X_1,\ldots,X_7)\;=\;c\cdot\det\left(
\frac{\partial^2 G}{\partial X_i\partial X_j}\right),
$$
(with $i,j\in\{1,2,3,6,7\}$) so the intersection $H_\alpha\cap \ccX$ is the 
Hessian ${\mathcal W}_{10}$ of Igusa's quartic.

The symmetric group $S_6\cong W(A_5)$ also acts on ${\mathcal W}_{10}$.
In fact, one has $b(\alpha,\alpha_i)=0$
for $i=2,\ldots,6$,
hence the root system of type $A_5$ defined by $\alpha_2,\ldots,\alpha_6$ 
is perpendicular to $\alpha$. 
Thus the Weyl group $W(A_5)$ acts on the hyperplane section 
${\mathcal W}_{10}=H_\alpha\cap \ccX$.

Now we find a 3-dimensional subdomain of $\ccS_4^M$ which maps to 
${\mathcal W}_{10}$.
As $W(E_6)$ acts transitively on the roots of $E_6$, we may redefine $\alpha:=\alpha_6$. 
The reflection $s_\alpha$ acts as $M_B\in Sp(8,\ZZ)$ 
(see Section \ref{actionMB}),
thus the fixed points of $M_B$ in $\ccS_4^M$ map to an eigenspace of 
$s_\alpha$ in $\PP^5$.

The fixed points in $\ccS_4$ of the involution $M_B$  are the $\tau\in \ccS_4$
such  that $M_B\cdot \tau=\tau$, that is, $B\tau B^{-1}=\tau$, equivalently, $B\tau=\tau B$,
so 
$$
\ccS_4^{M_B}\,=\,\left\{
\left(\begin{array}{cc}\tau_1 & \tau_2 \\ \tau_2&\tau_1 \end{array}\right)\;\in\ccS_4:\quad
\tau_1={}^t\tau_1,\quad \tau_2={}^t\tau_2\;\right\}~.
$$

The intersection of $\ccS_4^{M,M_B}:=\ccS_4^{M}\cap\ccS_4^{M_B}$
is three dimensional, because, as in Section \ref{weiltype}, we find that a period matrix in 
$\ccS_4^{M_B}$ is fixed
by $M$ iff $\tau_1=-(\tau_2+{}^t\!\tau_2)=-2\tau_2$. 
Thus, changing the sign of $\tau_2$, we get
$$
\ccS_4^{M,M_B}\,=\,\,\left\{
\left(\begin{array}{cc}2\tau_2 & -\tau_2 \\ -\tau_2&2\tau_2 \end{array}\right)\;\in\ccS_4:
\quad \tau_2\in \ccS_2\;\right\}~,
$$
in fact, $\tau_2$ must be symmetric and its imaginary part must be positive definite, conversely, given $\tau_2\in\ccS_2$ we get an element in $\ccS_4$
since the matrix with rows $2,-1$; $-1,2$
(the Cartan matrix of $A_2$) is positive definite.
For $\tau\in\ccS_4^{M,M_B}$, the abelian fourfold $X_\tau$ 
is isogeneous to 2 copies of the 
abelian surface $X_{\tau_2}$, since one has, similar to Section \ref{auto}, the following commutative diagram:
{\renewcommand{\arraystretch}{1.3}
$$
\begin{array}{ccccccccc}
0&\longrightarrow&\ZZ^4&\stackrel{\Omega_{\tau_2}}{\longrightarrow}&\CC^2
&\longrightarrow &X_{\tau_2}&\longrightarrow&0\\
&&N\downarrow\phantom{N}&&
N'\downarrow\phantom{N'}&&
\psi\downarrow\phantom{\psi}&&\\
0&\longrightarrow&\ZZ^8&\stackrel{\Omega_\tau}{\longrightarrow}&\CC^4
&\longrightarrow &X_\tau&\longrightarrow&0
\end{array}
\qquad\Bigl(\mbox{so}\quad N\Omega_{\tau}\,=\,\Omega_{\tau_2} N'\,\Bigr)~,
$$
}% end stretch 
where the maps are defined by matrices with $2\times 2$ blocks and the vectors are row vectors:
$$
N\,:=\,\,\left(\begin{array}{cccc}I& I&0&0\\0&0&I&I \end{array}\right),\quad
\Omega_{\tau}\,:=\,\left(\begin{array}{cc}
2\tau_2&-\tau_2\\
-\tau_2&2\tau_2\\
I&0\\0&I
\end{array}\right),\quad
\Omega_{\tau_2}\,:=\,\left(\begin{array}{cc}\tau_2\\ I\end{array}\right),\quad
N'\,:=\,\left(\begin{array}{cc}I&I\end{array}\right)~.
$$
This shows that there is a non-trivial holomorphic map $\psi:X_{\tau_2}\rightarrow X_\tau$.
Applying the automorphism $\phi$ of order three of $X_\tau$ 
to the image of $\psi$, 
one obtains another copy (up to isogeny) of the abelian surface $X_{\tau_2}$ 
in the fourfold $X_\tau$ and thus $X_\tau$ is isogeneous to $X_{\tau_2}^2$.
\qed

\subsection{Another reducible section of $\ccX$}\label{sing67}
In Proposition \ref{psing22} we showed that $\ccX\cap Z$, 
where $Z\cong \PP^3$ is subspace of $\PP^5$ perpendicular to 
a root subsystem of $E_6$ of type $A_2$, was reducible. 
Now we consider the intersection of $\ccX$ with a $\PP^3$ 
which is perpendicular to two orthogonal roots, 
so a root system of type $A_1^2$. 

Up to the action of $W(E_6)$ such a subsystem is unique, 
and there are $270$ of these.
We take the perpendicular roots $(0,0,0,0,1,1)$ and $(0,0,0,0,1,-1)$,
then the $\PP^3$ perpendicular to both roots is defined by $X_6=X_7=0$.

\subsection{Proposition}\label{psing67}
Let $W\subset \PP^5$ be the projective 3-space defined by
$$
W\,:\qquad X_6\,=X_7\,=\,0\qquad(\subset\,\PP^5)~.
$$
The intersection of $W$ with $\ccX\subset\PP^5$ has two irreducible components,
$\ccX\cap W=Q_{67}\cup S_{67}$, where $Q_{67}$ is a smooth quadric
and $S_{67}$ is a degree $8$ (singular)  rational surface.

The quadric $Q_{67}$ parametrizes abelian fourfolds whose endomorphism algebra contains the field of $12$-th roots of unity. The surface $S_{67}$
parametrizes abelian fourfolds which are isogeneous to a product of elliptic curves.
 
\ts
We found by explicit computation that
$$
\ccX\,\cap\,W\;=\;Q_{67}\,\cup\,S_{67},\qquad
Q_{67}\,:\;X_0^2 -X_1^2 - X_2^2 - X_3^2\,=\,0,
$$
and the degree $8$ surface $S_{67}$ is defined by
$$
X_0^2X_1^2X_2^2X_3^2 - X_1^4X_2^4 - X_1^4X_3^4  - X_2^4X_3^4 + X_1^4X_2^2X_3^2  + X_1^2X_2^4X_3^2 +     X_1^2X_2^2X_3^4\,=\,0~.
$$
The intersection of these two surfaces consists of $8$ conics. 
These conics are also the intersection of $Q_{67}$ with the planes 
$X_0=\pm X_1\pm X_2\pm X_3$, for any choice of signs.

The singular locus of $S_{67}$ consists of the three lines $X_i=X_j=0$ for 
$i,j\in\{1,2,3\}$. 
Let $S$ be the double cover  of $\PP^2$ (with coordinates $Z_0,Z_1,Z_2$)
given by  
$$
S:\quad
T^2=Z_0^4+Z_1^4+Z_2^4-Z_0^2Z_1^2-Z_0^2Z_2^2-Z_1^2Z_2^2
\qquad (S\subset W\PP(1,1,1,2))~.
$$
Then there is a birational isomorphism between $S_{67}$ and $S$ given by 
$$
S_{67}\, \dashrightarrow\, S,\qquad
(Z_0:Z_1:Z_2:T)\,:=\,(X_2X_3\,:\,X_1X_3\,:\,X_1X_2\,:\,X_0X_1X_2X_3)~,
$$
with birational inverse
$$
S\, \dashrightarrow\,  S_{67},\qquad
(X_0:X_1:X_2:X_3)\,:=\,(T\,:\,Z_1Z_2\,:\,Z_0Z_2\,:\,Z_0Z_1)~.
$$
Notice that the branch locus of $S\rightarrow \PP^2$ is reducible:
$$
Z_0^4+Z_1^4+Z_2^4-Z_0^2Z_1^2-Z_0^2Z_2^2-Z_1^2Z_2^2\,=\,
(Z_0^2 + \ww^2Z_1^2 + \ww Z_2^2)(Z_0^2 + \ww Z_1^2 + \ww^2 Z_2^2)~.
$$
The conics defined by the two factors intersect in the four points $(1:\pm 1:\pm 1)$. 
It is now easy to see that $S_{67}$ is rational: 
the inverse image of a general line through the point $(1:1:1)$ is again isomorphic 
to $\PP^1$ and each conic gives a (ramification) point on this double cover. 
Thus one can parametrize the double cover.

We now show that the quadric $Q_{67}$ 
parametrizes abelian fourfolds $X$ whose endomorphism algebra 
contains the field of $12$-roots of unity, $\QQ(\zeta_{12})\subset End_\QQ(X)$.

For this, it is convenient to change the pair of perpendicular roots to
$\alpha_3=(0,1,-1,0,0,0)$ and $\alpha_6=(0,0,0,0,1,-1)$. 
The projective 3-space $W'$ perpendicular to both of these roots is defined by $X_1=X_2$, $X_6=X_7$.
We write
$$
\ccX\,\cap\,W'\,=\,Q_{36}\,\cup\, S_{36}
$$ 
where $Q_{36}$ and $S_{36}$ 
are surfaces of degree $2$ and $8$ respectively. 

To find a subdomain of $\ccS_4^M$ which maps to $Q_{36}$,
we define an element $M_C\in Sp(8,\ZZ)$ as a block-matrix with four diagonal $2\times 2$ blocks $C={}^t\!C^{-1}$:
$$
C\,:=\,\left(\begin{array}{cc}0& 1\\-1&0 \end{array}\right)\quad(\in GL(2,\ZZ)),
\qquad
M_C\,:=\mbox{diag}(C,C,C,C)
\quad(\in Sp(8,\ZZ))~.
$$
One easily verifies that $M_CM=MM_C$, so $M_C$ lies in $C_M$, 
the centralizer of $M$. 
Thus $M_C$ maps $\ccS^M$ into itself. 
As $M$ has order $3$ and $M_C$ has order four (in fact $M_C^2=-I$ so $M_C$ has eigenvalues $\pm i$ with $i^2=-1$), 
the matrix $M_{12}:=MM_C$ is an element of order 
twelve in $Sp(8,\ZZ)$. 
As $M_{12}^4=M$, $M_{12}^3=M_C$, we get
$$
(\ccS^M)^{M_C}\,=\,\ccS^{M_{12}},\qquad M_{12}\,:=\,MM_C~.
$$
Using a commutative diagram as in Section \ref{auto}, one finds that for 
$\tau\in\ccS^{M_{12}}$ the matrix $M_{12}$ induces an automorphism $\phi_{12}$ of order $12$  the abelian variety $X_\tau$. 
Thus the field $\QQ(\zeta_{12})$, where $\zeta_{12}$ is a primitive $12$-root of unity, is contained in the endomorphism algebra of $X_\tau$ for any $\tau\in
\ccS^{M_{12}}$.
The eigenvalues of $\phi_{12}$ on the tangent space $T_0X_\tau=\CC^4$
are $\zeta_{12},\zeta_{12}^5,\zeta_{12}^7,\zeta_{12}^{11}$.

From \cite[Section 9.6]{BL} one then obtains 
(with $d=1$, $e_0=2$, $m=2$ and $(r_\nu,s_\nu)=(1,1)$ for $\nu=1,2$
and notice that their disc $\ccH_{1,1}$ is biholomorphic to $\ccS_1$) 
that abelian fourfolds with such an automorphism are parametrized by $\ccS_1^2$.
In particular, $\dim \ccS^{M_{12}}=2$. From \cite[Exercise 9.10 (4)]{BL} it
follows that for $\tau\in\ccS^{M_{12}}$ the endomorphism algebra of $X_\tau$  contains an indefinite quaternion algebra over the field $\QQ(\sqrt{3})$.

Now we show that the closure of the image of $\ccS^{M_{12}}
$ under $\Theta:\ccS_4^M\rightarrow\PP^5$ is the quadric $Q_{36}\subset Z'$.
The image is contained in the fixed points of $M_C$ acting on $\PP^{5}\subset\PP^{15}$. 
The action of $M_C$ on $\PP^{15}$ can be found as we did in Section \ref{actM} 
for $M$ (now with $A$ replaced by $\mbox{diag}(C,C)$),
one finds that the coordinates $X_0,X_3,X_{12},X_{15}$ are fixed
and that $M_C$ interchanges $X_1\leftrightarrow X_2$, $X_4\leftrightarrow X_8$,
$X_5\leftrightarrow X_{10}$, $X_6\leftrightarrow X_{9}$, $X_7\leftrightarrow X_{11}$, $X_{13}\leftrightarrow X_{14}$.
Thus:
$$
M_C:\,\PP^5\,\longrightarrow \,\PP^5,\qquad
(X_0:X_1:X_2:X_3:X_6:X_7)\,\longmapsto\,(X_0:X_2:X_1:X_3:X_7:X_6)~,
$$ 
hence $(\ccS^M)^{M_C}=\ccS^{M_{12}}$ maps to the subspace $W'=(\PP^5)^{M_C}\subset\PP^5$. 

The image of $\ccS^{M_{12}}$ is thus $Q_{36}$ or $S_{36}$. 
We checked that none of the quadrics $Q_m$'s is identically
zero on $W'$ and that there are exactly six which vanish on $Q_{36}$ (and as 
$S_{36}$ has degree $8$, none of the $Q_m$ can vanish on it).
Now we show that six of the $Q_m$ vanish on the image of $\ccS^{M_{12}}$,
hence $Q_{36}$ is the closure of the image of $\ccS^{M_{12}}$.

For this we use the action of $M_C$ on the $\theta[{}^{\epsilon}_{\epsilon'}]$,
similar to \cite[Proposition 10.7.3]{pmv}. 
Using the series defining these theta constants and the fact that $C{}^t\!C=I$ 
one has, with now $A=\mbox{diag}(C,C)$:
$$
\theta[{}^{\epsilon}_{\epsilon'}](M_C\cdot\tau)\,=\,
\theta[{}^{\epsilon A}_{\epsilon' A}](\tau)\,=\,
\theta\thchar{-\epsilon_2}{\epsilon_1}{-\epsilon_4}{\epsilon_3}
{-\epsilon_2'}{\epsilon_1'}{-\epsilon_4'}{\epsilon_3'}(\tau)\,=\,
(-1)^{\epsilon_2\epsilon_2'+\epsilon_4\epsilon_4'}
\theta\thchar{\epsilon_2}{\epsilon_1}{\epsilon_4}{\epsilon_3}
{\epsilon_2'}{\epsilon_1'}{\epsilon_4'}{\epsilon_3'}(\tau)~.
$$
It easily follows that if $M_C\cdot\tau=\tau$, then there are six $\theta_m(\tau)$ 
which satisfy $\theta_m(M_C\cdot\tau)=-\theta_m(\tau)$
and thus they vanish. 
Therefore six $Q_m$ vanish on the image of 
$\Theta(\ccS^{M_{12}})\subset Q_m$.

Finally we identify the period matrices mapping to the octic surface $S_{36}$. 
We define an element $M_D\in Sp(8,\ZZ)$ as a block-matrix with four diagonal $2\times 2$ blocks $D={}^t\!D^{-1}$:
$$
D\,:=\,\left(\begin{array}{cc}0& 1\\1&0 \end{array}\right)\quad(\in GL(2,\ZZ)),
\qquad
M_D\,:=\mbox{diag}(D,D,D,D)
\quad(\in Sp(8,\ZZ))~.
$$
One easily verifies that $M_DM=MM_D$, so $M_D$ lies in $C_M$, 
the centralizer of $M$. 
Thus $M_D$ maps $\ccS^M$ into itself. 
Moreover, $M_BM_D=M_DM_B$, thus $M_D$ maps   $\ccS_4^{M,M_B}$ 
into itself. 

The fixed point set $\ccS_4^{M,M_B,M_D}$ has dimension two because
it consists of the matrices $\tau=\tau(\tau_2)\in\ccS_4^{M,M_B} $ as in 
the proof of Proposition \ref{fixMB}, 
with $\mbox{diag}(D,D)\tau=\tau\mbox{diag}(D,D)$, so 
$$
\ccS_4^{M,M_B,M_D}\,=\,\left\{\left(\begin{array}{cc}2\tau_2 & -\tau_2 \\ -\tau_2&2\tau_2 \end{array}\right)\;\in\ccS_4:
\quad\tau_2\,=\,\left(\begin{array}{cc}\tau_{11}&\tau_{12}\\
\tau_{12}&\tau_{11}\end{array}\right)\;\in \ccS_2
\;\right\}~.
$$
In particular, $\dim \ccS_4^{M,M_B,M_D}=2$. As $M_D\equiv M_C$ mod $\Gamma_4(2,4)$, these matrices act in the same way on $\PP^{15}$ and thus 
$\Theta(\ccS_4^{M,M_B,M_D})$ is a surface in $W'$.

For $\tau\in \ccS_4^{M,M_B}$, the abelian variety $X_\tau$ 
is isogeneous to a selfproduct $Y_\tau^2$ where $Y_\tau$ is an abelian surface.  
The map $M_D$ induces an involution on the abelian variety $Y_\tau$, and thus
$Y_\tau$ is isogeneous to a product $E_1\times E_2$. Therefore if 
$\tau\in \ccS_4^{M,M_B,M_D}$, the abelian variety
$X_\tau$ is isogeneous to a product $(E_1\times E_2)^2$, where the $E_i$ are 
elliptic curves depending on $\tau$ and, for dimension reasons, 
the moduli of these two elliptic curves vary independently over $\ccS_1$. 
This implies that $End_\QQ(X_\tau)=M_2(\QQ)^2$ for a general $\tau\in\ccS_4^{M,M_B,M_D}$. The minimal polynomial of an element in this $\QQ$-algebra cannot be irreducible of degree $4$, and thus $\QQ(\zeta_{12})\not\subset
End_\QQ(X_\tau)$. Therefore the closure of $\Theta(\ccS_4^{M,M_B,M_D})$ 
is not $Q_{36}$, but it is the octic surface $S_{36}$.
\qed

\

\section{The boundary of $\ccX$}

\subsection{The boundary components}
The quotient of the Satake compactification $\overline{\ccA_4(2,4)}$ of
$\ccA_4(2,4)$
by the group $\Gamma_4:=Sp(8,\ZZ)$ is the Satake compactification $\overline{\ccA}_4=\overline{\ccS_4/\Gamma_4}$.
This variety has one boundary component of dimension $k(k+1)/2$ for $k=0,1\ldots,3$, whose closure is isomorphic to the Satake compactification $\overline{\ccA}_k:=\overline{\ccS_k/\Gamma_k}$.
Therefore the group $\Gamma_4$ acts transitively on the boundary components 
of a given dimension of $\overline{\ccA_4(2,4)}$.  
The boundary of $\ccX=
\Theta\big(\overline{\ccA_4(2,4)}\big)\cap\PP^5$, 
is, by definition, the intersection of $\PP^5$ with the boundary of 
$\Theta(\overline{\ccA_4(2,4)})$.

\subsection{Proposition}\label{bdprop}
\begin{enumerate}
\item
The boundary of $\ccX$ consists of $45$ lines, which we call the boundary lines.
They are the $W(E_6)$-orbit of the line 
$$
l\,:\quad X_2\,=\,X_3\,=\,X_6\,=\,X_7\,=\,0~.
$$

\item There are $27$ points of intersection of the boundary lines.
These points are called the cusps of $\ccX$. The group $W(E_6)$
acts transitively on the cusps.

\item The cusps  correspond to degenerate ppav's $(\CC^\times)^4$. 
A point on a boundary line which is not a cusp corresponds to a degenerate ppav $(\CC^\times)^2\times B$, where $B$ is an abelian surface with an automorphism of order three.

\item
Each boundary line contains $3$ cusps and each cusp is on $5$ of the boundary lines.
The $3$ cusps on $l$ are:
$$
(1:0:0:0:0:0),\quad
(1:1:0:0:0:0),\quad (1:-1:0:0:0:0)~.
$$
\end{enumerate}

\ts
We recall some facts on the action of $\Gamma_4$ on $\PP^{15}$. 
The normal subgroup $\Gamma_4(2)/\Gamma_4(2,4)\cong(\ZZ/2\ZZ)^8$ of 
$Sp(8,\ZZ)/\Gamma_4(2,4)$ is generated by block matrices
with $a,d=I$, and  $c=0$ and $b$ a diagonal matrix with even entries
or the transposed of such a matrix.
These matrices act as elements of a finite Heisenberg group on $\PP^{15}$. 
The matrix $M_{\beta,0}$ with $c=0$ and $b=\mbox{diag}(2\beta)$ 
acts as $X_\sigma\mapsto (-1)^{\beta{}^t\sigma}X_\sigma$
and $M_{0,\gamma},$ with $b=0$ and $c=\mbox{diag}(2\gamma)$,
acts as $X_\sigma\mapsto X_{\sigma+\gamma}$.

From  \cite[section 3.7]{pmv} 
(but note that we interchanged the diagonal blocks in $\tau(t)$)
it follows that a $k(k+1)/2$-dimensional boundary component of
$\Theta(\overline{\ccA_4(2,4)})$ is contained in the linear subspace
defined by $X_\sigma=0$ for those $\sigma\in(\ZZ/2\ZZ)^4$ with 
$(\sigma_{4-k},\ldots,\sigma_{4})\neq (0,\ldots,0)$.
This $\PP^{2^{k}-1}$ is a common eigenspace of the elements in the
Heisenberg group with $c=0$.
Thus to find all boundary components of $\ccX$ one determines the intersection of $\ccX$
with the eigenspaces of the elements of the Heisenberg group.

For $k=0$, one finds a zero dimensional boundary component of the image of 
$\overline{\ccA_4(2,4)}$ in $\PP^{15}$, it is the point $p:=(1:0:\ldots:0)$, 
so only $X_0\neq 0$. 
This point is fixed under the action of the subgroup  
in the Heisenberg group of matrices with $c=0$. 
The point $p$ actually lies in $\ccX\subset\PP^5$ and we checked that all 0-dimensional boundary components of $\ccX$ are in the $W(E_6)$-orbit of $p$, which has $27$ elements.

In terms of the $E_6$ root system, $p$ can be described as follows.
The fundamental weight $\lambda_i$ of $E_6$ is defined by the equations
$\bar{B}_0(\lambda_i,\alpha_j)=\delta_{ij}$, where $\delta_{ij}$ is Kronecker's delta.
It is easy to check that 
$
\lambda_2=(2,0,0,0,0,0),
$
thus $p$ is the image of this fundamental weight in $\PP^5=\PP(R(E_6)\otimes\CC)$.

The $\PP^3\subset\PP^{15}$ defined by 
by $X_{abcd}=0$ if $(c,d)\neq (0,0)\in\FF_2^2$ is the closure of a 3-dimensional boundary component of $\Theta(\overline{\ccA_4(2,4)})$.
Using the action of $\Gamma_4$, one finds that also the $\PP^3$ defined by
$X_{abcd}=0$ if $(a,c)\neq (0,0)\in\FF_2^2$ is the closure of a boundary component.
On this $\PP^3$, only the coordinates $X_0,X_1,X_4,X_5$ are non-zero 
and the other $12$ are zero. 
Using the results from Section \ref{actM}, one finds that this $\PP^3$ intersects 
the eigenspace $\PP^5$ in the line $l$ defined by $X_2=X_3=X_6=X_7=0$. 

By restricting the quadrics $Q_m$ to the line $l$, one finds that there are only three points in which $120$ of them vanish, thus there are exactly three cusps, 
the ones listed in the proposition, on $l$. 
The other results stated in the proposition follow from further Magma computations.
\qed

\

\

\end{document}